\documentclass[11pt,twoside]{article}

\title{From spatially periodic instantons to singular monopoles}   %
\date{}                                                            %
\author{Benoit Charbonneau
  \footnote{This research was supported by NSERC PGS A-B and MIT.}}%

\usepackage{fancyhdr}                                %
\usepackage{amssymb}                                 %
\usepackage{amsmath}                                 %
\usepackage{amscd}                                   %
\input diagrams  
\diagramstyle[size=1.8em,noPS]                       %
\usepackage{latexsym} 
\usepackage{pstcol,pst-node}                         %
\usepackage[hyperindex=true,                         %
pdfauthor={Benoit Charbonneau (benoit@alum.mit.edu)},%
  pdftitle={From spatially periodic instantons to    %
   singular monopoles}]{hyperref}                    %

\newlength{\myparskip}              \newlength{\myparindent}             %
\newlength{\oldparskip}             \newlength{\oldparindent}            %
\parskip\myparskip                  \parindent\myparindent               %
\fancypagestyle{plain}{%
   \fancyhead[LO]{Published in\\
              {\sc Communications in \\
               analysis and geometry}\\
               Volume 14, Number 1, 1-32, 2006} 
   \fancyhead[CO]{}\fancyhead[RO]{}
}

\newenvironment{proof}{\topsep=\smallskipamount \partopsep=0pt  %
 \begin{trivlist} \itemindent=\parindent                        %
 \item[\hskip \labelsep\emph{Proof:}]}{\qed\end{trivlist}}      %
\let\qed=\relax                                                 %
\def\qed                                                        %
 {{\unskip\nobreak\hfil\penalty50                               %
   \quad\hbox{}\nobreak\hfil $\Box$                             %
   \parfillskip=0pt \finalhyphendemerits=0 \par}}               %
\def\@thmcountersep{-}                                          %
\newtheorem{theorem}{Theorem}[section]                          %
\newtheorem{corollary}[theorem]{Corollary}                      %
\newtheorem{lemma}[theorem]{Lemma}                              %

\def\citep#1#2{\cite[{#1}]{#2}}

\DeclareMathAlphabet{\mathdj}{U}{msb}{m}{n}  

\DeclareMathAlphabet{\Gothique}{U}{euf}{m}{n}
\fontfamily{ptm}\selectfont

\DeclareMathSymbol{\sshortmid} {\mathord}{AMSb}{"70}
\DeclareMathSymbol{\sulcorner} {\mathord}{AMSa}{"70}
\DeclareMathSymbol{\surcorner} {\mathord}{AMSa}{"71}
\DeclareMathSymbol{\sllcorner} {\mathord}{AMSa}{"78}
\DeclareMathSymbol{\slrcorner} {\mathord}{AMSa}{"79}

\newcommand{\spinC}[1]{\,\,\,\makebox[0pt]{$#1$}
    \raisebox{1pt}{\makebox[0pt]{$\slash$}}\hspace{3pt}}   
\newcommand{\DD}{\spinC{\mathfrak{D}}}                     
\newcommand{\DDD}{{\spinC{\mathfrak{D}}^{!}}}     

\newcommand{\Z}{\ensuremath{\mathdj {Z}}} 
\newcommand{\R}{\ensuremath{\mathdj {R}}} 
\newcommand{\C}{\ensuremath{\mathdj {C}}} 

\newcommand{\TR}{{\R\times T^3}}  
\newcommand{\RT}{\TR}   
\newcommand{\Tp}{{[0,\infty)\times T^3}}
\newcommand{\Tm}{{(-\infty,0]\times T^3}}
\newcommand{\Rfs}{{\R^4}^*}   
\newcommand{\Rts}{{\R^3}^*}   
\newcommand{\Bw}{{B^3(w)}}  
\newcommand{\YHR}{[a,\infty)\times T^3}
\newcommand{\Ya}{\{a\}\times T^3}    
\newcommand{\Yap}{\YHR}        
\newcommand{\YR}{\R\times T^3}  

\newcommand{\SUn}{\mathrm{SU}(n)}        
\newcommand{\SUt}{\mathrm{SU}(2)}        

\newcommand{\coker}{\mathrm{coker}}
\renewcommand{\Im}{\mathrm{Im}}
\newcommand{\ind}{\mathrm{ind}}
\newcommand{\Spec}{{S\negthinspace p\negthinspace e\negthinspace c}}
\newcommand{\Vlam}{W}    
\newcommand{\Wlam}{{W'}}

\newcommand{\Ehat}{{V}}    
\newcommand{\Eh}{{V}}
\newcommand{\vE}{{V_{\negthickspace\negthickspace\lrcorner}}}
\newcommand{\bE}{{{}^{\raisebox{-1pt}{$\scriptstyle\ulcorner$}}\negthinspace\negthinspace V}}
\newcommand{\HHH}{{\mathcal{H}}} 
\newcommand{\Khat}{{K}}    
\newcommand{\vK}{{K_{\negthickspace\negthickspace
\lrcorner}}}
\newcommand{\bK}{{{}^{\raisebox{-1pt}{$\scriptstyle\ulcorner$}}\negthickspace\negthinspace K}}
\newcommand{\Ed}[1]{{V(#1)}}
\newcommand{\Kd}[1]{{K(#1)}}

\newcommand{\Ah}{{B}}   
\newcommand{\Higgs}{{\Phi}}
\newcommand{\Pb}{{\Phi^{\scriptscriptstyle\perp}}} 
\newcommand{\Pv}{{{\Phi}^{\slrcorner}}}           
\newcommand{\Prv}{{P^{\slrcorner}}}       
\newcommand{\Pp}{\Pb}                     
\newcommand{\BB}{{\overline{B}}}          
\newcommand{\bb}{{\Phi}}              
\newcommand{\bbb}{{\overline{B}}}  
\newcommand{\pp}{\eta}             

\newcommand{\Wotd}[1]{{W}^{1,2}_{#1}} 
\newcommand{\Wot}{{W}^{1,2}}          
\newcommand{\WH}{W^{\frac12,2}}       
\newcommand{\Ltl}{L^2_{\Lambda}}      
\newcommand{\Grid}{{\Gothique G}}
\newcommand{\iso}{\cong}              
\newcommand{\LHS}{left-hand-side}       
\newcommand{\RHS}{right-hand-side}      
\newcommand{\cqfd}{The proof is now complete.} 
\newcommand{\rk}{\mathrm{rk}}         
\newcommand{\End}{\mathrm{End}}       
\newcommand{\n}[2]{\|#1\|\strut_{#2}} 
\newcommand{\cl}{cl}                  
\newcommand{\scp}[1]{\langle #1\rangle} 
\newcommand{\grad}{{\scriptstyle g\negthinspace%
        r\negthinspace a\negthinspace d}\ }    
\newcommand{\dhat}{{d^z}}             
\newcommand{\limit}{\partial\text{-term}}
\newcommand{\bigdist}[2]{{\mathrm{dist}\bigl(#1,#2\bigr)}} 

\newcommand{\ulepsi}{{\left\ulcorner\!\!\epsilon\right.}}
\newcommand{\uepsi}{{\overline{\epsilon}}}
\newcommand{\urepsi}{{\left.\epsilon\!\!\right\urcorner}}
\newcommand{\lepsi}{{\sshortmid\epsilon}}
\newcommand{\repsi}{{\epsilon\sshortmid}}
\newcommand{\eb}{{\ulepsi}}
\newcommand{\ev}{{{}_\llcorner\negthickspace\epsilon}}
\newcommand{\ep}{{\uepsi}}
\newcommand{\emo}{{\underline\epsilon}}
\newcommand{\elr}{{\epsilon_{\negthickspace\lrcorner}}}
\newcommand{\eur}{{\urepsi}}
\newcommand{\eul}{\eb}
\newcommand{\el}{{\lepsi}}
\newcommand{\er}{{\repsi}}

\newcommand{\dll}{{{}_\sllcorner\negthickspace\delta}}
\newcommand{\dlr}{{\delta_{\negthinspace\negthickspace\slrcorner}}}
\newcommand{\dur}{{\delta\negthickspace\raisebox{3pt}{$\scriptstyle\surcorner$}}}
\newcommand{\dul}{{{}^{\raisebox{-2pt}{$\scriptstyle\sulcorner$}}\negthickspace\delta}}

\begin{document}
\thispagestyle{plain}
\maketitle


\begin{abstract}
The main result is a computation of the Nahm transform of a $\SUt$-instanton
over $\R\times T^3$, called spatially-periodic instanton. It is a singular 
monopole over $T^3$, a solution to the Bogomolny equation, whose rank is 
computed and behavior at the singular points is described.
\end{abstract}

\section{Introduction}               
Heuristically, there is a correspondence, called the \emph{Nahm
transform}, between
\begin{enumerate}
\item solutions to the anti-self-dual (ASD) equation, or its the
appropriate dimensional reduction, on the quotient of $\R^4$ by a
closed subgroup $\Lambda$ of $\R^4$, and satisfying a finite energy
condition, and
 \item solutions to some associate equation satisfying some boundary
condition on the quotient of ${\R^4}^*$ by the dual subgroup
$\Lambda^*=\{ f\in {\R^4}^*\mid f(\Lambda)\subset\Z\}$.
\end{enumerate}

This heuristic comes from a re-engineering due to Nahm \cite{nahm}
of the ADHM construction of instantons on $\R^4$ \cite{ADHM}.
Nahm's approach has the advantage of being transportable to
quotients by non trivial subgroup $\Lambda$ as well,  with some ad
hoc efforts necessary in each case.

Nahm gave an outline of the correspondence for classical
instantons ($\Lambda=\{0\}$) and for monopoles on $\R^3$ ($\Lambda
= \R$). Corrigan--Goddard in \cite{corrigan-goddard} completed the
details of the ADHM construction following Nahm's guideline, while
Hitchin in \cite{hitchinMonopoles}  completed the story for
$\SUt$-monopole on $\R^3$.  In \cite{nakajima}, Nakajima rendered
Hitchin's proof more parallel to the ADHM story.

This framework guided several other authors in the quest for an
understanding of other moduli spaces of instantons (or their
appropriate dimensional reduction) on various quotients of $\R^4$:
for instantons on $T^4$, see \cite{schenk1988,braam1989}; for
monopoles for other classical groups, see \cite{hurtubisemurray};
for calorons, or instantons on $S^1\times \R^3$, see
\cite{singernye,nye}; for instantons on $T^2\times \R^2$, see
\cite{jardim1999,jardim2001,jardim2002b,jardim2002,biquardjardim};
and for monopoles on $\R^2\times S^1$, see
\cite{cherkis1998,cherkis1999,cherkis2001}.  Marcos Jardim wrote a
survey paper \cite{jardimsurvey} on the Nahm transform,
and the reader is invited to consult it for some insights on an
even more general framework in which to place the above referenced
literature and the present paper.

Apart from some numerical approximations and remarks in
\cite{vanBaal1996} and a computation of the Nahm transform of charge $1$ instantons in
\cite{vanBaal1999}, the case of the \emph{spatially periodic
instantons}, instantons on $\RT$, has been largely ignored.  The
present paper starts the groundwork necessary to close that gap.
We prove here that the Nahm transform of an instanton on $\RT$ is
a singular monopole on $T^3$ with specific behavior at the
singular points.

This paper is organized as follows.  The main result on the Nahm
transform of instantons on $\R\times T^3$ and its singular
behavior is  spelled out in Section \ref{sec:result} after the
adequate language is explained.  Before reaching this result, it
is useful to go over a brief overview of the classical ADHM
construction in Section \ref{sec:ADHMN}, then check the bigger
picture of the Nahm transform heuristic in Section
\ref{chap:heuristic}, and then zoom in on the Fredholmness
properties of the Dirac operators on $\R\times T^3$ in Section
\ref{sec:Fredholm}.  The proof of the result splits three ways:
first, the rank of the transformed bundle is computed at the end
of Section \ref{sec:result}; then, a splitting of the transformed
bundle around the singularities is developed in Section
\ref{sec:exactseq}; and finally, the asymptotic of the Higgs field
is proved in Section \ref{sec:Higgsfield}.

\emph{Acknowledgments: This work is part of my Ph.D.
thesis \cite{benoitthesis}.  Grateful for the guidance and support of
my advisor Tomasz Mrowka, I thank him warmly. I also thank
Larry Guth, Peter Kronheimer, Fr\'ed\'eric Rochon, and
Michael Singer for stimulating discussions.}

\section{The classical ADHMN}                  
\label{sec:ADHMN}                          
The classical work of Atiyah, Drinfeld, Hitchin and Manin,
termed ADHM construction, classify
all the solutions to the ASD equation on $\R^4$, up to gauge
equivalence.  Once viewed under the umbrella of the Nahm transform
heuristic, thus adding an N to form \emph{ADHMN}, the classification
is as follows.

A connection $A$ on a $\SUn$-bundle $E$ over $\R^4$ whose
curvature $F_A$ satisfies the ASD equation $*F_A=-F_A$ and the
finite energy condition $\int_{\R^4}|F_A|^2<\infty$ gives rise, through an analysis of its Dirac operator $\DD_A$,  to
a set of algebraic data: two vector spaces
\begin{gather*}
V=L^2\cap \ker(\DD^*_A),\text{ and}\\
W=\text{ bounded harmonic sections of $E$ for
$\nabla_A$},\end{gather*} and five maps
\begin{gather*}
\bb_1,\ldots,\bb_4\colon V\to W,\\
\pp\colon V\to S^+\otimes W. \end{gather*}

Since the vector space is built using the augmented Dirac operator
$\DD_A^*$ acting on sections of $S^-\otimes E$, the dimension of
$V$ can be computed by some index theorem, and
\[\dim V =\frac{1}{8\pi^2}\int_{\R^4}|F_A|^2\]
provided the cokernel $L^2\cap \ker(\DD_A)$ is $\{0\}$.  It is
indeed so, as the Weitzenbock formula
\[\DD_A^*\DD_A=\nabla_A^*\nabla_A+\cl(F_A^+)\]
clearly establishes:  for an instanton connection, the Clifford
multiplication term vanishes and a $L^2$ solution $\phi$ to
$\DD_A^*\phi=0$ must be parallel, hence $0$ since $\R^4$ has
infinite volume.

The map $\bb_i=Pm_{x_i}$ is the composite of the multiplication by
the $i$th coordinate, denoted $m_{x_i}$ and the $L^2$-projection
$P$ on $\ker(\DD_A^*)$, while the map $\pp$ encodes  the
asymptotic behavior of elements of $V$.

For an instanton $(E,A)$, the associated algebraic data
$(V,W,\bb,\pp)$ satisfy a non-degeneracy condition and the ADHM
equation, the precise formulation of which is not important here.
This ``ADHM transform" places in one-to-one
correspondence instantons modulo gauge equivalence with
non-degenerate solutions to the ADHM equation modulo some symmetry
group action. A complete description of this construction can be
found in \cite[Chap. 3]{DK}, and in the author's thesis
\cite{benoitthesis}.

It is a fruitful idea to interpret the set of maps
$\bb=(\bb_1,\ldots,\bb_4)$ as a constant connection form
\[\bbb=\bb_1 dx^1+\cdots \bb_4 dx^4\] on the trivial bundle
$\underline V$ over $\R^4$ with fiber $V$.  The curvature $F_\bbb$
of $\bbb$ splits as
\[F_\bbb = (\text{ASD part}) + (\text{ SD part involving }\pp).\]

Morally, the idea is that the transformed connection $\bbb$ on
${\R^4}^*$, invariant under the action of $\Lambda^*={\R^4}^*$, is
almost anti-self-dual, and the self-dual part is determined by the
asymptotic behavior of harmonic spinors.

\section{The Nahm transform heuristic}                
                 \label{chap:heuristic}        
The work of Nahm provides a framework in which to think about the
classification of all the finite energy solutions to the ASD
equation on a quotient $\R^4/\Lambda$. Philosophically, once we
find the appropriate codomain for the Nahm transform to be
described in this section, it should be an isomorphism. This idea
has been shown to work in many cases, as explained in the
introduction.

A connection $A$ on a $\SUn$-bundle $E$ over $\R^4$, invariant
under the action of a closed subgroup $\Lambda$, and whose
curvature $F_A$ satisfies the ASD equation \[*F_A=-F_A\] and the
finite energy condition \[\int_{\R^4/\Lambda}|F_A|^2<\infty\]
gives rise, this time, to a bundle $V$ with a connection $B$ over
$\Rfs/\Lambda^*$. Those objects are constructed in the following
way.

For  an element $z$ of $\Rfs$, the space of $\R$-valued linear
functions on $\R^4$, we define the bundle $L_z$ \index{$L_z$} over
$\R^4$ to be a trivial $\C$-bundle with connection
\[\omega_z:=2\pi i z= 2\pi i\sum_{j=1}^4 z_j dx^j.\]
For $z'\in\Lambda^*$, the flat bundles $L_z$ and $L_{z+z'}$ over
$\R^4/\Lambda$, both invariant under the action of $\Lambda$, are
isomorphic. We write $A_z$ for the connection $A\otimes 1+1\otimes
\omega_z$ on $E\otimes L_z=E$.  For $z\in\Rfs$, consider the
operator
\[\DD_{A_z}^*\colon \Gamma(\R^4,S^-\otimes E\otimes L_z)\to
   \Gamma(\R^4,S^+\otimes E\otimes L_z).\]
A section of the bundle $S^-\otimes E\otimes L_z$
is said to be in $L^2_\Lambda$ if it is
invariant under the action of $\Lambda$ and if its $L^2$-norm over
$\R^4/\Lambda$ is finite.

The first ingredient of the Nahm transform of the instanton
$(E,A)$ is the family of vector spaces
\[\Eh_z:=L^2_\Lambda\cap\ker(\DD_{A_z}^*).\]

Since the vector space $V_z$ is built using the augmented Dirac
operator $\DD_{A_z}^*$ acting on sections of $S^-\otimes E$, the
dimension of $V_z$ can often be computed by an appropriately
chosen index theorem, and it is  constant on connected components
on which $\DD_{A_z}$ is Fredholm provided the cokernel
$L^2_\Lambda\cap \ker(\DD_{A_z})$ is $\{0\}$.  For a quotient
$\R^4/\Lambda$ of infinite volume, it is indeed so, as the
Weitzenbock formula
\[\DD_A^*\DD_A=\nabla_A^*\nabla_A+\cl(F_A^+)\]
clearly establishes:  for an instanton connection, the Clifford
multiplication term vanishes and a $L^2$ solution $\phi$ to
$\DD^*_A\phi=0$ must be parallel, hence $0$ because of the
infinite volume condition.  For a quotient of finite volume, we
must add an extra condition to ensure  the cokernel is trivial.

It turns out in many cases that $\DD_{A_z}^*$ is not Fredholm for
every $z$, which is a good thing.  Suppose for example that $\DD_{A_z}^*$
was Fredholm everywhere
when $\Lambda=\Z^3$. As we explore in this present paper, the
object created by the Nahm transform is a monopole over
$T^3$.  But as one can show (see \cite[Prop. 1]{pauly}), smooth monopoles
over compact $3$-manifolds\index{monopoles} are not very
interesting.

Set $g_z(x):=e^{2\pi i z(x)}$. 
Notice that for any section $\phi$ of $S^-\otimes E$, we have
$\DD_{A_z}^*(g_z\phi)=g_z\DD^*_{A} \phi$. Then for all
$z'\in\Lambda^*$, we have an isomorphism
\begin{equation}\label{isomgauge}
   g_{z'}\colon \Eh_z\to \Eh_{z+z'},\end{equation}
hence $\Eh$ is a bundle over $\Rfs/\Lambda^*$.

\subsection{First viewpoint: on $\Rfs$, a curvature computation}
In the understanding of the ADHM construction, it was beneficial
to  view the maps $\bb_i$ as parts of a connection on the bundle
$\underline{\Eh}$ on $\Rfs$, without passing to the quotient.  We
do similarly here and consider first the bundle $\Eh$ on an open
subset of $\Rfs$ on which the Dirac operator is Fredholm.

We define  a connection $\Ah$ on $\Eh$.  Each fiber $\Eh_z$ is in
fact contained in the vector space $\Ltl(S^-\otimes E)$.  We can then consider the
trivial connection $\dhat$ in the trivial bundle of fiber
$\Ltl(S^-\otimes E)$, and its projection $P\dhat$ to $\Eh$.

The operator $\DD^*_{A_z}\DD_{A_z}$ should be invertible, and we
use its inverse, the Green's operator
$G_{A_z}=(\DD^*_{A_z}\DD_{A_z})^{-1}$, to define the projection $P$
by the formula
\[P=1-\DD_{A_z}G_{A_z}\DD_{A_z}^*.\]

To parallel the ADHMN story, let's now compute the curvature $F_\Ah$ of $\Ah$.  
To simplify
the notation, we set $\Omega:=2\pi i\sum_{j=1}^4\cl(dx^{j})dz^{j}.$
Then $[\dhat,\DD_{A_z}]=\Omega$, and similarly for $\DD_{A_z}^*$.

The curvature $F_\Ah$ can be computed as follows:
\begin{align*}
    \scp{(P\dhat)^{2}\phi,\psi} &=\scp{\dhat P\dhat\phi,\psi}\\
         &=\scp{P\dhat\phi,\dhat\psi}-\scp{\dhat\phi,\dhat\psi}\\
       &=-\scp{\DD_{A_z}G_{A_z}\DD_{A_z}\dhat\phi,\dhat\psi}\\
        &=\scp{\DD_{A_z}G_{A_z}\Omega\phi,\dhat\psi}.
\end{align*}

Let $\nu$ be the normal vector field to $S^{r-1}(R)\times T^s$.
The integration by parts necessary to bring $D$ on the \RHS\ of
the scalar product introduces a boundary term
\begin{equation}\label{eqn:limit}
\limit:=\lim_{R\to \infty}\int_{S^{r-1}(R)\times T^s}
\scp{\cl(\nu)G_{A_z}\Omega\phi,\dhat\psi}.\end{equation}

Performing the said integration by parts, we obtain
\begin{align*}\scp{F_{\Ah} \phi,\psi}
        &=\scp{G_{A_z}\Omega\phi,\DD_{A_z}\dhat\psi}+\limit\\
        &=-\scp{G_{A_z}\Omega\phi,\Omega\psi}+\limit\\
        &=\scp{G_{A_z}\phi,\Omega\wedge\Omega\psi}+\limit.
\end{align*}

In terms of the usual basis $\epsilon_j$  and $\bar{\epsilon}_j$
of  respectively $\bigwedge^+$ and $\bigwedge^-$, we have
\[\Omega\wedge \Omega=-4\pi^{2}\sum_{j=1}^3\bigl(\cl(\epsilon_{j})\epsilon_{j}+
    \cl(\bar\epsilon_{j})\bar\epsilon_{j}\bigr).\]
Since $\bigwedge^{+}$ acts trivially on $S^{-}\otimes E$, the
first term of the curvature is ASD.

In the case we are studying at this moment, the $\limit$  is
$0$.

\subsection{Second viewpoint: on a $3$-dimensional quotient, the
Bogolmolny equation} Let's now shift our perspective and look at
$V$ and $\BB$ from the viewpoint of the quotient.  Suppose some
$\R$ is in $\Lambda^*$, say as the axis $z_1$. In fact, suppose
here that $\Lambda=\Z^3$, and thus that $\Lambda^*=\R\times\Z^3 $.
Set $g_z(x)=e^{2\pi i x_1z_1}$.  Then
\[g(\BB)=-2\pi i P m_{x_1} dz^1+P\bigl( \frac\partial{\partial z_2} dz^2+\cdots +
\frac\partial{\partial z_4} dz^4
\bigr).\]

So using this gauge transformation, we render $\BB$ independent of
the $z_4$ coordinate.  We define the Higgs field $\Higgs$ by
\[\Higgs = -2\pi i Pm_{x_1},\]
and the connection $B$ on $\Rfs/\Lambda^*=T^3$ by
\[\Ah=Pd^z,\]
where $z$ represents here the coordinates $(z_2,z_3,z_4)$ on
$T^3$.  As we just saw,
\[g(\BB)=\Higgs dz^1+\Ah.\]
Should we be able to prove that $\limit=0$, it would be so that
$g(\BB)$ is ASD.  It is in fact so, as we see in the next section,
and thus $(\Ah,\Higgs)$ satisfies the dimensional reduction of the
ASD equation
\[\nabla_\Ah\Higgs = *F_\Ah\]
called the \emph{Bogomolny equation}.

\section{Fredholmness of the Dirac operator}             
\label{sec:Fredholm}                              
It is crucial now to understand exactly for which $z\in T^3$ the Dirac
operator $\DD_{A_z}^*$ acting on $L^2$ sections of $S^-\otimes E$ over
$\RT$ is Fredholm.

Let's start with a $\SUt$-instanton $(E,A)$ on $\RT$ and call $t$
the $\R$-coordinate.  Modulo gauge transformation, we can pick a
representative in temporal gauge: $A$ has no $dt$ term and can be
seen as a path of connections on $T^3$, parameterized by $\R$.  In
temporal gauge, the Dirac operator splits as
\[\DD_{A}^*=-\frac\partial{\partial t} +D_{A}\]
with $D_A$ the Dirac operator on the cross-section $\{t\}\times T^3$.
Furthermore, as $t\to\infty$ and $t\to-\infty$, the connection $A$ has
flat limits $\Gamma_+$ and $\Gamma_-$; see \cite[Thm 4.3.1]{mrowkabook}.  Consequently, the operator
$D_{A_z}$ limits to $D_{\Gamma_{+\ z}}$ and $D_{\Gamma_{-\ z}}$ at
$+\infty$ and $-\infty$.  It is a crucial observation of
Atiyah--Patodi--Singer \cite{APS1} that the unbounded operator
$\DD_{A_z}^*\colon L^2\to L^2$ is Fredholm if and only if $0$ is not
in the spectrum of either $D_{\Gamma_{+\ z}}$ or $D_{\Gamma_{-\ z}}$;
see \cite[Chap. 6]{benoitthesis} for a very detailed account.

As it turns out, any flat $\SUt$ bundle over a $3$-manifold splits as
a sum of flat $U(1)$-bundles.  Our bundle $E$, restricted to
$\pm\infty$, splits respectively as
\[E=L_{w_\pm}\oplus L_{-w_\pm},\]
for some $w_\pm\in \Rts$.  The spectrum of $D_{\Gamma_{+\ z}}$ is thus
the multiset
\begin{equation}\Spec(D_{\Gamma_{+\, z}})=\pm2\pi\bigl|\Lambda^*_\Z-w_+-z\bigr|
\cup\pm2\pi\bigl|\Lambda^*_\Z+w_+-z\bigr|\label{spec}\end{equation}
for the part $\Lambda^*_\Z\iso\Z^3$ of $\Lambda^*$ in $\Rts$, and similarly
for $D_{\Gamma_{-\,z}}$; see \cite[Chap. 3]{benoitthesis}.

Thus, $\DD_{A_z}^*$ is Fredholm as long as $z$ is not in the set
\[W=\{w_+,-w_+,w_-,-w_-\}.\]
Keeping a parallel with the notation for the ADHMN story, the set $W$ is in
some sense our set of ``infinity data,''  although in a much milder
way than for $\R^4$.

It is appropriate at this point to ask for which $z$ is $\DD_{A_z}^*$
Fredholm when we change the domain to allow for more or less
growth.  Choosing a weight $\delta\in\R^2$, say
$\delta=(\delta_-,\delta_+)$, and a weighing function $\sigma_\delta$
such that
\[\sigma_\delta=\begin{cases} e^{-\delta_- t}, &\text{ for } t<-1,\\
                  e^{-\delta_+ t}, &\text{ for } t>1,\end{cases}\]
we define the weighted $L^2$-norm
\[\n{f}{L^2_\delta}:=\n{\sigma_\delta f}{L^2},\]
and naturally
\[L^2_\delta:=\bigl\{ f\in L^2_{loc}\mid \n{f}{L^2_\delta}<\infty\bigr\}.\]
We omit the bundle from the notation, as it should always be clear which bundle is involved.

Similarly, we can define weighted Sobolev spaces.  These include only
those $L^2_\delta$ sections whose derivatives are also in $L^2_\delta$.
Fix a 
connection $\nabla$ on $E$, and set
\[\Wotd\delta:=\{f\in L^2_\delta\mid \nabla f\in L^2_\delta\}.\]

Keeping in mind that the first coordinate of the weight describes the
growth at $-\infty$ while the second describes the growth at $+\infty$, we
define the grid
\[\Grid_A:=\Spec(D_{\Gamma_-})\times\R\ \cup \ \R\times\Spec(D_{\Gamma_+})\]
in the weight space $\R^2$.  Naturally, the Atiyah--Patodi--Singer condition
becomes
\[\DD^*_{A_z}\colon \Wotd\delta\to L^2_\delta\
\text{ is Fredholm if and only if }\ \delta\not\in\Grid_{A_z}.\]

We define the spaces
\begin{equation}\begin{aligned}
  \ker(\delta)&:=\ker(\DD_A\colon \Wotd\delta \to L^2_\delta),\\
  \ker^*(\delta)&:=\ker(\DD^*_A\colon \Wotd\delta \to L^2_\delta),
\end{aligned}\label{def:ker}\end{equation}
and the integers
\begin{equation}\begin{aligned}
    \ind(\delta)&:=\ind (\DD_A\colon \Wotd\delta \to L^2_\delta)\\
     N(\delta)&:=\dim\ker(\delta),\text{ and}\\
     N^*(\delta)&:=\dim\ker^*(\delta).
\end{aligned}\label{def:N}\end{equation}

Since $(L^2_\delta)^* = L_{-\delta}^2$, elliptic regularity tells us
that $\dim\coker(\DD_A) = N^*(-\delta)$, hence
\[\ind(\delta)=N(\delta)-N^*(-\delta).\]
That the formal adjoint $\DD^*_A$ on $\Wotd{-\delta}$ is really the adjoint
of $\DD_A$ on $\Wotd\delta$ is guaranteed by the following lemma.

\begin{lemma}\label{lemma:kill}
The subspace $\ker^*(-\delta)$ of $L^2_{-\delta}={(L^2_\delta)}^*$ kills
$\Im(\delta)$ in the $L^2$ natural pairing.
\end{lemma}

\begin{proof}
Suppose $\phi$ is a smooth function with compact support.  Then for all
$\psi\in\ker^*(-\delta)$, we have $\scp{\psi,\DD\phi}=\scp{\DD^*\psi,\phi}=0$.
Since $C_c^\infty$ is dense in $\Wotd\delta$, the lemma holds.
\end{proof}

The operator $\DD_{A_z}^*\colon \Wotd\delta\to L^2_\delta$ is
conjugate to the operator $\DD_{A_z}^*+\sigma_\delta\cl(\grad
\sigma_\delta^{-1})$ {}from $\Wot$ to $L^2$.  So the family
parameterized by $\delta$ in an open square delimited by
$\Grid_{A_z}$  is continuous and hence has constant index.  In
fact, the dimensions of the kernel and the cokernel are also
constant in an open square.  The proof is easy and can be found in 
\cite[Thm 6.3-2]{benoitthesis}.

As we cross a wall in $\Grid_A$ to change from one open square to another,
the index $\ind$ of $\DD_A$ and the index $\ind^*$ of $\DD_A^*$ change as follows:
\begin{align}
   \ind(\delta)&=\ind(\eta)+\dim\{D_{\Gamma_+}\phi=-\lambda\phi\},
        \text{ and}\nonumber\\
   \ind^*(\delta)&=\ind^*(\eta)+\dim\{D_{\Gamma_+}\phi=
        \lambda\phi\}\label{wallcrossingeqn}
\end{align}
when $\delta_+<\eta_+$, and $\delta$ and $\eta$ are in adjacent open
squares separated by the wall $\R\times\{\lambda\}\subset\Grid_A$;
\begin{align*}
\ind(\delta)&=\ind(\eta)+\dim\{D_{\Gamma_-}\phi=-\lambda\phi\},\text{ and}\\
\ind^*(\delta)&=\ind^*(\eta)+\dim\{D_{\Gamma_-}\phi=\lambda\phi\}
\end{align*}
when $\delta_->\eta_-$, and $\delta$ and $\eta$ are in adjacent open
squares separated by the wall $\{\lambda\}\times\R\subset\Grid_A$.

When the limit $\Gamma_+$ is such that the kernel of $D_{\Gamma_+}$ is
$\{0\}$,
not only are $\DD_A$ and $\DD_A^*$ Fredholm, as we saw above, 
we also have that
$A$ decays exponentially to $\Gamma_+$.  So there exist $\beta>0$ such that
$|A-\Gamma_+|\leq C e^{-\beta t}$ for $t>0$;  this is a consequence
of \cite[Thm 5.2.2]{mrowkabook} and of the embedding of
$\Wot$ in bounded $C^0$ functions, \cite[Thm 3.4]{hebey}.  In that case, we have the
following result on harmonic spinors.

\begin{theorem}\label{thm:asympt}               
Suppose $\phi\in \ker(\DD^*_A)\cap \Wotd\delta$.
Suppose $\lambda-\beta<\eta<\delta$ and that $\lambda$ is the only
eigenvalue of $D_{\Gamma_+}$ between $\eta$ and
$\delta$: $\Spec(D_{\Gamma_+})\cap[\eta,\delta]=\{\lambda\}$.    Then there
exist an eigenvector $\bar\psi$ of $D_{\Gamma_+}$ of eigenvalue $\lambda$
on $T^3$ and $\bar\phi\in\Wotd\eta((0,\infty)\times T^3)$ such that
\begin{equation}\label{eqn:asympt}
   \phi=e^{\lambda t}\bar \psi +\bar\phi\text{ for }t>0.
\end{equation}
Furthermore, $\bar\phi=O(e^{\eta t})$ as $t\to\infty$.
\end{theorem}

\begin{proof}
The space $L^2(T^3)$ splits according to the finite dimensional eigenspaces
$\Vlam_\lambda$ for $D_{\Gamma_+}$.  Let $\Pi_\delta^+$, $\Pi_\delta^-$ and
$\Pi_\delta$ be respectively
the projections from $L^2(T^3)$ to 
\[\bigoplus_{\lambda>\delta}\Vlam_\lambda,\ \ 
\bigoplus_{\lambda<\delta}\Vlam_\lambda, \text{ and }
\Vlam_\delta.\]
To simplify notation
we omit $\delta$ when it is $0$ and 
set $\phi^\pm:=\Pi^\pm(\phi)$.

For every $\phi\in L^2(Y)$, let $\phi_\lambda$ be its $\Vlam_\lambda$ 
component.
Thus $\phi=\sum \phi_\lambda.$
Using this decomposition, we can define the space $\WH(T^3)$ using the norm
\begin{equation}\label{def:normhalf}\n{\phi}{\WH}^2=\sum(1+|\lambda|)\n{\phi_\lambda}{L^2}^2.\end{equation}

Because $T^3$ is compact, the space $\WH(T^3)$ defined by two different
Dirac operators are equal, with commensurate norms.  The $+$ and $-$ part of
$L^2$, however, depend highly on $D_{\Gamma_+}$.

The operator
\begin{equation}
\begin{aligned}
\DDD_{\Gamma_+}\colon \Wot(\Yap)&\to L^2(\Yap)\oplus\Pi^+\WH(\Ya)\\
  \phi&\mapsto (\DD_{\Gamma_+}\phi,\Pi^+\phi(a))
\end{aligned}\label{Dhalf}
\end{equation}
is an isomorphism when $D_{\Gamma_+}$ has no kernel.

The proof of this fact starts as one does in the full cylinder case:
\begin{align*}
  \n{\DD_{\Gamma_+}\phi}{L^2}^2&=\n{\partial_t\phi}{L^2}^2
              +\n{D_{\Gamma_+}\phi}{L^2}^2
              +\int_a^\infty\partial_t\scp{\phi,D_{\Gamma_+}\phi}_{L^2(Y)}\\
&\geq C\n{\phi}{\Wot}^2-\scp{\phi(a),D_{\Gamma_+}\phi(a)}_{L^2(Y)}.
\end{align*}
Contrary to the full cylinder case, the boundary term here cannot be made to vanish and henceforth helps
control the $\Wot$-norm of $\phi$.  Using the inequality above and 
the decomposition $\phi=\sum \phi_\lambda$, we find
\begin{equation}
\n{\phi}{\Wot}^2 
\leq C\bigl(\n{\DD_{\Gamma_+}\phi}{L^2}^2+\n{\phi^+(a)}{\WH(T^3)}^2\bigr).\label{ineq:ellEstHalf}
\end{equation}
We just proved that $\n{\phi}{\Wot}\leq C\n{\DDD_{\Gamma_+}\phi}{}$, 
hence $\DDD_{\Gamma_+}$ is 
semi-Fredholm and injective.  Suppose now  that $(\psi,\eta)$ is 
perpendicular to $\Im(\DDD_{\Gamma_+})$.  
For all  $\phi\in \Wot(\YHR)$, we have 
\begin{align*}
0&=\scp{\DD_{\Gamma_+}\phi,\psi}+\scp{\eta,\phi^+(a)}\\
 &=\scp{\phi,\DD_{\Gamma_+}^*\psi}-\scp{\phi(a),\psi(a)}+\scp{\eta,\phi^+(a)}\\
 &= \scp{\phi,\DD_{\Gamma_+}^*\psi}-\scp{\phi^-(a),\psi^-(a)}+\scp{\eta-\psi^+(a),\phi^+(a)}.
\end{align*}
Going  through all the $\phi$ with $\phi(a)=0$ in a first time, $\phi^+(a)=0$ then, and finally $\phi^-(a)=0$, we 
prove
\begin{gather*}\DD_{\Gamma_+}^*\psi=0,\\ \eta=\psi^+(a),\\ \psi^-(a)=0.\end{gather*}

Thus we have $-\partial_t\psi+D_{\Gamma_+}\psi=0$, which means that $\psi$ is a linear 
combination of the $e^{\lambda t}\psi_\lambda$.
The condition $\psi^-(a)=0$ forces out all the negative $\lambda$, while the
positive ones are forced out by the $L^2$ condition.  Hence $\psi=0$ and
$\DDD_{\Gamma_+}$ is surjective.  The proof that the operator (\ref{Dhalf})
is an isomorphism is now complete.  

For $a$ big enough, the  operator $\DDD_{A}$,
not independant of $t$ but close enough to $\DDD_{\Gamma_+}$, 
is also an
isomorphism.


As in the full cylinder case, we can look at weighted version of $\DD$ and
$\DDD$.  For computing the asymptotic expansion of
harmonic spinors, we actually need to consider the dual $\DD^*$ and its counterpart 
\begin{align*}
\DDD_A^*\colon \Wotd\delta(\Yap)&\to L^2_\delta(\Yap)\oplus \Pi^-_\delta\WH(T^3) \\
   \phi&\mapsto (\DD_A^*\phi,\Pi^-\phi(a)),\end{align*}
which is Fredholm if and only if $\delta\not\in \Spec(D)$, and is an 
isomorphism when Fredholm.

We close the proof of Theorem \ref{thm:asympt} with a diagram chase. 
We first introduce maps to compose our diagram.

Recall $\eta<\delta$ and $\Spec(D)\cap[\eta,\delta]=\{\lambda\}$.  
Then obviously,
the map 
\begin{align*}
   I\colon  \Pi_{\eta}^-\WH(\Ya)\oplus \Vlam_\lambda &\to \Pi_{\delta}^-\WH(\Ya)\\
             (\phi,\psi)&\mapsto \phi+e^{a\lambda}\psi
\end{align*}
is an isomorphism, while the map
\begin{align*}
   J\colon  \Wotd\eta(\Yap)\oplus \Vlam_\lambda &\to \Wotd\delta(\Yap)\\
             (\phi,\psi)&\mapsto \phi+e^{\lambda t}\psi
\end{align*}
is an injection.

Consider now the map
\begin{align*}
   K\colon  \Wotd\eta(\Yap)\oplus \Vlam_\lambda 
       &\to L^2_\eta(\Yap)\oplus \Pi_{\eta}^-\WH(\Ya)\oplus \Vlam_\lambda\\
             (\phi,\psi)&\mapsto 
                    \bigl(\DD_A(\phi+e^{\lambda t}\psi),\Pi_\eta^-\phi,
                 \psi+e^{-a\lambda}\Pi_\lambda\phi(a)\bigr).
\end{align*}
As $\bigl|\DD^*_A(e^{\lambda t}\psi)\bigr|\leq Ce^{(\lambda-\beta)t}|\psi|$, 
then $\DD^*_A(e^{\lambda t}\psi)\in L^2_\eta(\Yap)$, and $K$ is well-defined.

We put all these maps in a commutative diagram
\begin{equation}\label{diag:DIJK}
\begin{CD} \Wotd\delta(\Yap) @>{\DDD_A^*}>> 
          L^2_\delta(\Yap)\oplus \Pi_\delta^-\WH(\Ya)\\
@A{J}AA            @AA{\iota\oplus I}A\\
\Wotd\eta(\Yap)\oplus \Vlam_\lambda @>>K> L^2_\eta(\Yap)\oplus \Pi_\eta^-\WH(\Ya)\oplus \Vlam_\lambda
\end{CD} 
\end{equation}

We know that $\DDD_A^*$ is an isomorphism.  Using the identification 
\[\DDD_A^*\colon \Wotd\eta(\Yap)\equiv L^2_\eta(\Yap)\oplus \Pi_\eta^-\WH(\Ya),\] we see
that $K$ has the form
\[\begin{bmatrix}1&p\\q&1\end{bmatrix}\]
for the splitting $\Wotd\eta(\Yap)\oplus \Vlam_\lambda$ of the domain and codomain.
Hence $K-1$ is a compact operator, and $K$ is thus Fredholm of index $0$.
If $K(x)=K(y)$, then $\DDD_A^* J(x)=\DDD_A^* J(y)$ as the diagram is commutative,
hence $x=y$ and $K$ is injective.  Being of index $0$, it henceforth must
be an isomorphism.

Let's now exploit this fantastic diagram.  Suppose
\[\phi\in \ker(\DD^*_A)\cap \Wotd\delta(\YR).\]
Then  for $a$ big enough, the diagram (\ref{diag:DIJK}) has rows which 
are isomorphism for $\delta$ and $\eta$ satisfying the hypothesis 
of the theorem.

We now chase around the diagram.  Since $I$ is an isomorphism, we know there
exist $(\chi,\nu)\in \Pi_\eta^-\WH(\Ya)\oplus \Vlam_\lambda$ such that
\[\iota\oplus I(0,\chi,\nu)=\DDD_A^*(\phi).\]
But as $K$ is an isomorphism, there is 
$(\bar\phi,\bar\psi)\in \Wotd\eta(\Yap)\oplus \Vlam_\lambda$ such that
\[K(\bar\phi,\bar\psi)=(0,\chi,\nu).\]
By commutativity of the diagram, we have
\[\DDD_A^* J(\bar\phi,\bar\psi) =\DDD_A^*(\phi)\]
but $\DDD_A^*$ is an isomorphism hence
   $\phi=e^{\lambda t}\bar \psi +\bar\phi$ for $t>a$.

As the choice of $a$ is artificial, we set $a=0$.  
The proof is now complete.\end{proof}

Suppose now
\begin{gather*}
    \lambda\in \Spec(D_{\Gamma_-})\times \Spec(D_{\Gamma_+}),\\
    \delta\text{ is in the upper left open square adjacent to }\lambda,\\
    \eta\text{ is in the lower right open square adjacent to }\lambda.
\end{gather*}

When $A$ decays exponentially to its limits, we have
\begin{equation}\label{eqn:remark}\ker(\lambda)=\ker(\eta).\end{equation}
Indeed, suppose now $\phi\in\ker(\lambda)$.  Then $\phi\in\ker(\delta)$ hence
by Theorem \ref{thm:asympt}, we expand $\phi$ for $t>0$ as
$\phi=e^{-\lambda_+ t}\psi_{\lambda_+}
+\bar\phi$, with $\bar\phi\in\Wotd{\eta_+}([0,\infty)\times T^3)$.
Since $\phi$ and $\bar\phi$ are  both in $\Wotd{\lambda_+}$, so is the
term $e^{-\lambda_+ t}\psi_{\lambda_+}$.  This fact implies that
$\psi_{\lambda_+}=0$.  Using a similar proof at $-\infty$, we
find $\phi\in\Wotd\eta$.
Obviously, the same is true for $\ker^*$.

\section{Nahm Transform: Instantons to singular monopoles} 
\label{sec:result}                               
  Since
$\DD_{A_z}^*$ is Fredholm $L^2\to L^2$ outside of $W$, and
since $\ker(\DD_{A_z})=0$ as $F_{A_z}$ is
ASD and $\RT$ has infinite volume, we have a bundle $\Ehat$
over $T^3\setminus W$ whose fiber at $z$ is
\[\Ehat_z:=\ker(\DD_{A_z}^*)\cap L^2.\]

As outlined earlier, this bundle is equipped with
\begin{align*}
 &\text{a connection }\Ah \text{ on } T^3\setminus W,\\
 &\text{a Higgs field }\Higgs\in\Gamma(T^3\setminus W,\End \Ehat).
\end{align*}

The main result of this present paper is the following theorem.

\begin{theorem}\label{thm:SingularMonopole}
Outside of a set $W$ consisting of at most four points, the family of
vector spaces $\Ehat$ described above defines a vector bundle of rank
\[\frac1{8\pi^2}\int|F_A|^2,\]
and the couple
$(\Ah,\Higgs)$ satisfies the Bogomolny equation
\[\nabla_{\Ah}\Higgs=*F_\Ah.\]

For $w\in W$ and  $z$ close enough to  $w$,
there are maps $\Pb$ and $\Pv$ such that
\[\Higgs = \frac{-i}{2|z-w|}\Pb+\Pv,\]
and $\Pb$ is the $L^2$-orthogonal projection on the
orthogonal complement of a naturally defined subbundle $\vE$ of $\Ehat$.
\end{theorem}

The last part of the theorem is made clearer by the introduction of
$\vE$ in Section \ref{sec:exactseq}.

\begin{proof}
The rank of $\Ehat$ is computed in Lemma \ref{thm:LtwoIndex} below.

The boundary term of Equation (\ref{eqn:limit}) is
\[\limit=\scp{\nu\Omega G\phi,\dhat\psi}_{T^3}\Bigl|^\infty_{-\infty}.\]
For $z\not\in W$, both $G\phi$ and $\dhat\psi$ decay exponentially by
Theorem \ref{thm:asympt} hence
\[\limit=0,\]
and the connection $P\dhat$ on
$\R\times (T^3\setminus W)$ is ASD.
Thus,
as explained in Section \ref{chap:heuristic}, the pair $(\Ah, \Higgs)$
satisfies outside of $W$ the appropriate dimensional reduction of
the ASD equation, which is
in this case the Bogomolny Equation
\[\nabla_{\Ah}\Higgs=*F_{\Ah}.\]

The last part of the theorem is the content of Section \ref{sec:Higgsfield}
and rests on the  splitting of Section
\ref{sec:exactseq}.
\end{proof}

As announced, we compute now the rank of $V$, and prove an $L^2$-index theorem for $\R\times T^3$.

\begin{lemma}\label{thm:LtwoIndex}               
For a $SU(2)$-instanton $(E,A)$ on $\TR$, the index of the Dirac
operator
\[\DD_A\colon \Wot(\TR)\to L^2(\TR)\]
when $A$ has nonzero limits at $\pm\infty$ is given by the formula
\[\ind(\DD_A)=-\frac1{8\pi^2}\int|F_A|^2.\]
\end{lemma}

\begin{proof}                              
The fact that $A$ has nonzero limits guarantees that the
operator $\DD_A$ is Fredholm on $\Wot$.  Moreover, $A$ decays exponentially
to its limits.

Let \[(\chi_R^+,\chi_R^-,\chi_R^0)\]
be a partition of unity subordinate to the covering
\[\Bigl((R,\infty)\times T^3,  (-\infty,-R)\times T^3,
       (-R-1,R+1)\times T^3\Bigr).\]
Suppose $\Gamma_\pm =d +\gamma_\pm$, and $A=d+a$.  Then $a$
tends to $\gamma_+$ and $\gamma_-$ when
$t$ tends to $+\infty$ and $-\infty$ respectively.    Set
\begin{equation}\label{eqn:defaR}
  a_R=\chi_R^+\gamma_++\chi_R^-\gamma_-+\chi_R^0a. \end{equation}

The sequence $\DD_{a_n}-\DD_{a_R}$ of compact operators is 
Cauchy, thus has a limit, $K$ say, which is then compact.  
As $\DD_A=\DD_{a_R}+K$, we 
have that 
$\ind(\DD_A)=\ind(\DD_{a_R})$ for all $R>0$.  We now compute
$\ind(\DD_{a_R})$ using the relative index theorem.  It could be that
$\Gamma_-\neq\Gamma_+$, but this case is easily converted to a situation where
$\Gamma_-=\Gamma_+$, as we now see.

Choose a path $\Gamma_s$ in the space of flat connections on
$T^3$ starting at $\Gamma_+$ and ending at $\Gamma_-$, and avoiding the trivial
connection.  Hence $0\not\in \Spec(D_{\Gamma_s})$ for all $s$; recall Equation
(\ref{spec}).
Suppose $\Gamma_s=d +\gamma_s$ and set
\begin{equation}\label{eqn:defasR}
  a^s_R=\chi_R^+\gamma_s+\chi_R^-\gamma_-+\chi_R^0a. \end{equation}

The family $\DD_{a_R^s}$ of Fredholm operator depends continuously on $s$.
Hence
\[\ind(\DD_A)=\ind(\DD_{a_R})=\ind(\DD_{a_R^0})=\ind(\DD_{a_R^1}).\]
Note that the connection $a_R^1$ equals $\Gamma_-$ outside
$[-R-1,R+1]\times T^3$.  Hence the relative index theorem  tells us
\begin{equation}\label{eqn:relative}
    \ind(\DD_{a_R^1})-\ind(\DD_{\Gamma_-})=\ind(\tilde\DD_{a_R^1})-
     \ind(\tilde\DD_{\Gamma_-}),\end{equation}
where the tilded operators are extensions to some compact manifold of
the restriction of the operators $\DD_{a_R^1}$ and $\DD_{\Gamma_-}$
to $[-R-1,R+1]\times T^3$.

Because $D_{\Gamma_-}$ has no kernel, 
$\DD_{\Gamma_-}\colon \Wot\to L^2$ is an isomorphism, and thus
$\ind(\DD_{\Gamma_-})=0$.  Hence the \LHS{}
of Equation (\ref{eqn:relative}) is equal to $\ind(\DD_A)$.

To compute the \RHS, we embed $[-R-1,R+1]\times T^3$ in some flat $T^4$.
The spinor bundles $S^+$ and $S^-$ on $[-R-1,R+1]\times T^3$ agree very
nicely with those of  $T^4$.  We extend both $a_R^1$ and $\Gamma_-$ by
the trivial bundle with connection $\Gamma_-$.

The Atiyah--Singer index theorem
tells us that
\begin{align*}
\ind(\tilde\DD_{\Gamma_-})
        &=\bigl\{ch(\Gamma_-)\cdot\mathbf{\hat A}(T^4)\bigr\}[T^4]\\
\ind(\tilde\DD_{a_R^1})
        &=\bigl\{ch(a_R^1)\cdot\mathbf{\hat A}(T^4)\bigr\}[T^4]\\
        &=\bigl(\frac{c_1^2}2-c_2\bigr)[T^4].
\end{align*}

Since $a_R^1$ is in $\SUt$, we have $c_1=0$, while
\[c_2[T^4]=\frac1{8\pi^2}\int_{T^4}\bigl(|F_{a_R^1}^-|^2-|F_{a_R^1}^+|^2\bigr).\]

Note that on the complement of $[-R-1,R+1]\times T^3$ in $T^4$, the connection
$a_R^1$ equals $\Gamma_-$ hence is flat there.  Furthermore,
on $[-R,R] \times T^3$, we have $a_R^1=A$.  On $[R,R+1]\times T^3$ and
$[-R-1,-R]\times T^3$, the curvature $F_{a_R^1}$ involves cut off functions,
their derivatives and $(A-\Gamma_-)$ terms.  Since $A$ tends to $\Gamma_-$
exponentially fast, we therefore have constant $C$ and $\beta$ such that
\[\Bigl|\ind(\DD_A)+\frac1{8\pi^2}\int_{[-R,R]\times T^3}|F_A|^2\Bigr| \leq
C e^{-\beta R}.\]
As $R\to\infty$, we have the wanted result.
\end{proof}

\section{A Geometric Splitting and Exact Sequences}          
\label{sec:exactseq}                                   
In this section, we analyze a splitting of $\Ehat$ in a
neighborhood of a point $w\in W$ where the solution $(\Ah,\phi)$
to Bogomolny equation is singular.  This point $w$ is associated,
say, to the limit $\Gamma=\Gamma_+$ of $A$ at $+\infty$, in the
sense that $\Gamma$ splits $E$ as $L_w\oplus L_{-w}$ on $T^3$.

Suppose the connection $A$ decays at most with rate $\beta$, as in
$|A-\Gamma_+|\leq Ce^{-\beta t}$ for $t>0$ and
$|A-\Gamma_-|\leq Ce^{\beta t}$ for $t<0$.
Set
\[\epsilon:=\frac14\min
   \Bigl(\beta,\bigdist{w}{\Lambda^*+W\setminus\{w\}}\Bigr),\]
and
define the six weights \index{weights}
\begin{alignat*}6
\eb &:=(-\epsilon,\epsilon)\quad\quad &\ep &:=(0,\epsilon)
    \quad\quad &\eur&:=(\epsilon,\epsilon)\\
\ev &:=(-\epsilon,-\epsilon) &\emo &:=(0,-\epsilon)
    \quad\quad &\elr&:=(\epsilon,-\epsilon)
\end{alignat*}
displayed here in a way which is reminiscent of their position in $\R^2$.

Consider the ball $\Bw$ of radius $2\epsilon$ around $w$.  As $z$ varies in
$\Bw$, and depending on whether $\Gamma_+=\Gamma_-$ or not, there are two or
one walls to cross to pass from $0$ to $\eul$ and from $\elr$ to $0$.
In a picture, we have
\[\begin{pspicture}(-1.5,-2)(7.7,2)
\psline(-1.2,-0.5)(1.2,-0.5)\psline(-1.2,0.5)(1.2,0.5) 
\rput(2.1,0.55){$2\pi|z-w|$}\rput(2.1,-0.45){$-2\pi|z-w|$} 
\psline(-0.5,-1.1)(-0.5,1.2)\psline(0.5,-1.2)(0.5,1.2) 
\rput(1,1.4){$2\pi|z-w|$}\rput(-0.9,1.4){$-2\pi|z-w|$} 
\rput(-1,1){$\eul$}\rput(0,0){$0$}\rput(1,-1){$\elr$}
\rput(.2,-1.88){$\Gamma_+=\Gamma_-$}
\rput(4.8,0){%
\psline(-1.4,-0.5)(1.4,-0.5)\psline(-1.4,0.5)(1.4,0.5) 
\rput(2.4,0.55){$2\pi|z-w|$}\rput(2.4,-0.45){$-2\pi|z-w|$} 
\psline(-1.2,-1.2)(-1.2,1.2)\psline(1.2,-1.2)(1.2,1.2) 
\rput(-1,1){$\eul$}\rput(0,0){$0$}\rput(1,-1){$\elr$}
\rput(.2,-1.88){$\Gamma_+\neq\Gamma_-$}
}\end{pspicture}\]

As $z$ varies in $\Bw$, those walls move around without ever touching $\elr$
and $\eul$.  Hence for $L^2_\elr$ and $L^2_\eul$, the
operators $\DD_{A_z},\DD^*_{A_z}$
and $\DD^*_{A_z}\DD_{A_z}$ are Fredholm for all $z\in\Bw$.

Hence for $z\in \Bw$, the
six vector spaces
\begin{alignat*}5
\bE_z&:=\ker(\DD_{A_z}^*)\cap L^2_{\eul},
    \quad\quad&\bK_z&:=\ker(\DD_{A_z})\cap L^2_{\eul} ,
   \quad \quad\\
\vE_z&:=\ker(\DD_{A_z}^*)\cap L^2_{\elr} ,
    &\Khat_z&:=\ker(\DD_{A_z})\cap L^2\\
\HHH_z&:=\ker(\nabla_{A_z}^*\nabla_{A_z})\cap L^2_{\eul},\quad\quad &\vK_z&:=\ker(\DD_{A_z})\cap L^2_{\elr},
\end{alignat*}
are kernels of Fredholm operators.
By contrast, the space $\Ehat_z$, already defined as
$\ker(\DD_{A_z}^*)\cap L^2$, is not the kernel of a Fredholm operator at $w$.

Notice that none of those vector spaces form a priori a bundle over $\Bw$
as the dimensions could jump at random.  However, for
$L^2_{\eul}$
and $L^2_{\elr}$, the operators  $\DD_{A_z}$, $\DD_{A_z}^*$,
and
$\nabla_{A_z}^*\nabla_{A_z}$ are Fredholm operators for all $z\in \Bw$.  The
various indices are therefore constant and we have that, for example,
\[\dim \vE_z-\dim\bK_z\text{ is constant on }\Bw.\]

We have the following obvious results:
\begin{alignat*}2
&\vE\subset\Ehat\subset\bE,\quad\quad\quad & \vK\subset\Khat\subset\bK,\\
&\DD\HHH\subset\bE, &  \bK\subset \HHH,\\
&&\vK=\Khat=\{0\}.
\end{alignat*}

Equation (\ref{eqn:remark}) signifies here that $\vE_w=\Ehat_w$.
The following few lemmas describe in more detail the relationship
between the various spaces.

The smallest eigenvalues of
$D_{\Gamma_z}$ are $\pm 2\pi |z-w|$.  For simplicity, we set
\[\lambda:=2\pi|z-w|,\]
and define
\[\Vlam_\lambda:=\lambda\text{ eigenspace of }D_{\Gamma_z}\text{ on }T^3.\]

The family $\Vlam_\lambda$ defines a bundle over the sphere $|z-w|=\lambda/2\pi$
around $w$.  Its rank is given by
\begin{equation}\label{rkV}
\rk \Vlam_\lambda=\begin{cases}
   1,& \text{ if }\lambda\neq 0\text{ and }2w\not\in\Lambda^*;\\
   2,& \text{ if }\lambda\neq 0\text{ and }2w\in\Lambda^*,
       \text{ or }\lambda=0\text{ and }2w\not\in\Lambda^*;\\
   4,& \text{ if }\lambda= 0\text{ and }2w\in\Lambda^*.
\end{cases}\end{equation}
This $\Vlam_\lambda$ plays an important
role in understanding the relations between the various spaces just
introduced.

For any instanton connection $A'$ on $\TR$, set
\begin{align*}
   \Ed\delta &:=\ker(\DD^*_{A'})\cap L^2_\delta,\\
   \Kd\delta &:=\ker(\DD_{A'})\cap L^2_\delta,
\end{align*}
and let $[\delta]$ denote the open square in $\R^2\setminus\Grid_{A'}$
containing $\delta$.

\begin{lemma}[one wall]\label{lemma:onewall}Suppose $\delta,\eta\in\R^2\setminus\Grid_{A'}$
are weights for which $[\delta]$ and $[\eta]$ are adjacent and
separated by the wall $\{\mu\}\times \R$ or $\R\times\{\mu\}$.  Then the
sequence
\begin{equation}\label{seq:atw}\begin{CD}
  0\longrightarrow\Ed\delta\longrightarrow\Ed\eta
  @>{\scriptscriptstyle\lim (e^{-\mu t}\cdot)}>>\Vlam_\mu @>{\scriptscriptstyle\bigl(\lim (e^{\mu t}\cdot)\bigr)}^*>>
  \Kd{-\delta}^*\longrightarrow\Kd{-\eta}^*\longrightarrow0,
\end{CD}\end{equation}
where the limits are both evaluated at $+\infty$ when $[\eta]$ is
above $[\delta]$ and at $-\infty$ when $[\eta]$ is to the left of $[\delta]$,
is exact.
\end{lemma}

\begin{proof}Theorem \ref{thm:asympt} ensures that
the limits give functions $\alpha$ and $\beta^*$ which
are well defined, and that
\[0\longrightarrow\Ed\delta\longrightarrow\Ed\eta\longrightarrow \Vlam_\mu\quad
\text{ and }\quad
0\longrightarrow\Kd{-\eta}\longrightarrow\Kd{-\delta}\longrightarrow \Vlam_\mu\]
are exact.

It only remains to prove that Sequence (\ref{seq:atw}) is exact
at $\Vlam_\mu$.  Suppose $\phi\in \Ed\eta$ and $\psi\in \Kd{-\delta}$.
Then
\begin{align*} 0&= \scp{\DD_{A'}^*\phi,\psi}-\scp{\phi,\DD_{A'}\psi}\\
  &=\lim_{t\to\infty}\scp{\phi,\cl(\partial_t)\psi}-\lim_{t\to-\infty}\scp{\phi,\cl(\partial_t)\psi}\\
  &=\lim_{t\to\infty}\scp{e^{-\mu t}\phi,\cl(\partial_t) e^{\mu t}\psi}
    -\lim_{t\to-\infty}\scp{e^{-\mu t}\phi,\cl(\partial_t) e^{\mu t}\psi}.\end{align*}

One of those limits is $\beta^*\alpha(\phi)(\psi)$ while the other one vanishes
as we now see.  Suppose $[\eta]$ is above $[\delta]$, and suppose
$\{\mu'\}\times \R$ is the wall to their right.  Then $\phi=O(e^{\mu' t})$
as $t\to-\infty$ by Theorem \ref{thm:asympt}.  But for some $\mu''<\mu'$,
the wall $\{-\mu''\}\times \R$ is exactly to the right of $[-\eta]$ hence
$\psi=O(e^{-\mu'' t})$ as $t\to -\infty$.  But then
\[\beta^*\alpha(\phi)(\psi)=\lim_{t\to-\infty}O(e^{(\mu'-\mu'')t})=0,\]
hence $\Im(\alpha)\subset\ker(\beta^*)$.  A similar argument establish the
same fact when $[\eta]$ is to the left of $[\delta]$.

The sequence is then exact if
$\dim\Im(\alpha)=\dim\ker(\beta^*)$.
We have two short exact sequences:
\begin{gather*}0\longrightarrow\Ed\delta\longrightarrow\Ed\eta\longrightarrow
 \Im(\alpha)\longrightarrow0,\quad \text{ and }\\
 0\longrightarrow \Vlam_\mu/\ker(\beta^*)\longrightarrow\Kd{-\delta}^*
 \longrightarrow \Kd{-\eta}^*\longrightarrow0.\end{gather*}
Using those short exact sequences and notation
{}from Equations (\ref{def:N}), we have
\begin{align*}
  \dim\Im(\alpha)-\dim\ker(\beta^*)
   &=N^*(\eta)-N^*(\delta)-\dim \Vlam_\mu+N(-\delta)-N(-\eta)\\
   &=\ind^*(\eta)-\ind^*(\delta)-\dim \Vlam_\mu.\end{align*}
The Wall Crossing Equation
(\ref{wallcrossingeqn}) forces the last line to be $0$.
The proof is thus complete.\end{proof}

\begin{corollary}
Suppose $\Gamma_+\neq \Gamma_-$. Then the sequences
\begin{align}
0\longrightarrow \Ehat_z\longrightarrow \bE_z\longrightarrow \Vlam_\lambda
   \longrightarrow 0, &\quad\quad\text{ for }\lambda\neq0,\label{seq:hatbarz}\\
0\longrightarrow \vE_z\longrightarrow \Ehat_z\longrightarrow \Vlam_{-\lambda}
  \longrightarrow \bK_z\longrightarrow 0,
       &\quad\quad\text{ for }\lambda\neq0,\label{seq:hatV}\\
0\longrightarrow \Ehat_w \longrightarrow \bE_w\longrightarrow \Vlam_0
  \longrightarrow \bK_w\longrightarrow 0,& \label{seq:hatbarw}
\end{align}
are exact.
\end{corollary}

\begin{proof} Apply Lemma \ref{lemma:onewall} to the        
choice of weights
$\{\eul,0\}$ and $\{0,\elr\}$ for the connection $A'=A_z$, and remember
that $\vK=\Khat=\{0\}$.
\end{proof}

\begin{corollary}\label{lemma:seqtwowalls}           
Suppose $\Gamma_+= \Gamma_-$. Then the sequences
\begin{align}
0\longrightarrow \Ehat_z\longrightarrow \bE_z\longrightarrow \Vlam_\lambda\oplus \Vlam_{-\lambda}
   \longrightarrow 0, &\quad\quad\text{ for }\lambda\neq0,\label{seq2:hatbarz}\\
0\longrightarrow \vE_z\longrightarrow \Ehat_z\longrightarrow \Vlam_\lambda\oplus \Vlam_{-\lambda}
  \longrightarrow \bK_z\longrightarrow 0,
       &\quad\quad\text{ for }\lambda\neq0,\label{seq2:hatV}\\
0\longrightarrow \Ehat_w \longrightarrow \bE_w\longrightarrow \Vlam_0\oplus \Vlam_0
  \longrightarrow \bK_w\longrightarrow 0,& \label{seq2:hatbarw}
\end{align}
are exact.
\end{corollary}

\begin{proof}                                       
Suppose we have the following choice of weights:
\[\begin{pspicture}(-1.2,-1.2)(1.2,1.2)
\psline(-.9,0)(.9,0) 
\rput(1.1,0){$\mu$} 
\psline(0,-.9)(0,.9)
\rput(0,-1.1){$-\mu$}
\rput(.4,-.4){$\dlr$}\rput(.4,.4){$\dur$}
\rput(-.4,-.4){$\dll$}\rput(-.4,.4){$\dul$}
\end{pspicture}\]

Denote $\iota$ any inclusion map, and $L_\mu^\pm$ the maps
\[L_\mu^+(\phi)=\lim_{t\to\infty} e^{\mu t}\phi,\quad\text{ and } \quad L_\mu^-=\lim_{t\to-\infty} e^{\mu t}\phi.\]
Then sequences akin to Sequence (\ref{seq:atw}) fit in a  
diagram
\begin{equation}\label{tresse}
\begin{diagram} 
0&&\rTo&&\Ed\dur&&\rTo^{L_{-\mu}^+}&&\Vlam_\mu&&\rTo^{{L_\mu^+}^*}&&\Kd{-\dll}^*&&\rTo&&0\\ 
&\rdTo&&\ruTo_{\iota}&&\rdTo^{\iota}&&\ruTo_{L_{-\mu}^+}&&\rdTo^{{L_\mu^+}^*}&&\ruTo_{\iota}&&\rdTo^{\iota}&&\ruTo&\\
& &\Ed\dlr&&    &&\Ed\dul&&&&\Kd{-\dlr}^*&&&&\Kd{-\dul}^*&&\\
&\ruTo&&\rdTo_{\iota}&&\ruTo^{\iota}&&\rdTo_{L_\mu^-}&&\ruTo^{{L_{-\mu}^-}^*}&&\rdTo_{\iota}&&\ruTo^{\iota}&&\rdTo&\\
0&     &\rTo&     &\Ed\dll&&\rTo_{L^-_\mu}&&\Vlam_{-\mu}&&\rTo_{{L_{-\mu}^-}^*}&&\Kd{-\dur}^*&&\rTo&&0
\end{diagram}.\end{equation}

Suppose $\phi\in \Ed{\dul}$, and $\psi\in \Kd{-\dlr}$.  Then
\begin{align*} 0&=\scp{\DD_{A'}\phi,\psi}-\scp{\phi,\DD_{A'}\psi}\\
  &= \scp{\phi,\cl(\partial_t)\psi}|^{\infty}_{-\infty}\\
  &=\lim_{t\to \infty} \scp{e^{-\mu t}\phi,\cl(\partial_t) e^{\mu t}\psi}
   -\lim_{t\to -\infty} \scp{e^{\mu t}\phi,\cl(\partial_t) e^{-\mu t}\psi}\\
 &=\bigl({L_\mu^+}^*L_{-\mu}^+(\phi)-{L_{-\mu}^-}^*L_{\mu}^-(\phi)\bigr)(\psi),
\end{align*}
hence the middle square commutes.  It is quite obvious that all the
other squares and triangles commute.  From Diagram (\ref{tresse}), we extract,
for an obvious choice of maps, 
the exact sequence
\[0\longrightarrow \Ed\dlr \longrightarrow \Ed\dul
   \longrightarrow \Vlam_\mu\oplus \Vlam_{-\mu} \longrightarrow \Kd{-\dlr}^*
   \longrightarrow \Kd{-\dul}^*\longrightarrow 0.\]

In particular, the sets of weights
\[\begin{pspicture}(-1.5,-2)(8.6,2)
\psline(-1.2,-0.5)(1.2,-0.5)\psline(-1.2,0.5)(1.2,0.5) 
\rput(1.5,0.55){$\lambda$}\rput(1.5,-0.45){$-\lambda$} 
\psline(-0.5,-1.1)(-0.5,1.2)\psline(0.5,-1.2)(0.5,1.2) 
\rput(0.5,1.4){$\lambda$}\rput(-0.5,1.4){$-\lambda$} 
\rput(-1,1){$\eul$}\rput(0,1){$\ep$}
\rput(-1,0){$\el$}\rput(0,0){$0$}
\rput(0,-1.5){at $z\neq w$}
\rput(.2,-1.88){$\Gamma_+=\Gamma_-$}
\rput(4,0){%
\psline(-1.2,-0.5)(1.2,-0.5)\psline(-1.2,0.5)(1.2,0.5) 
\rput(1.5,0.55){$\lambda$}\rput(1.5,-0.45){$-\lambda$} 
\psline(-0.5,-1.1)(-0.5,1.2)\psline(0.5,-1.2)(0.5,1.2) 
\rput(0.5,1.4){$\lambda$}\rput(-0.5,1.4){$-\lambda$} 
\rput(1,0){$\er$}\rput(0,0){$0$}
\rput(1,-1){$\elr$}\rput(0,-1){$\emo$}
\rput(0,-1.5){at $z\neq w$}
\rput(.2,-1.88){$\Gamma_+=\Gamma_-$}
}
\rput(7.8,0){%
\psline(-1,0)(1,0) 
\rput(1.2,0){$\lambda$} 
\psline(0,-1)(0,1)
\rput(0,1.2){$-\lambda$}
\rput(.5,-.5){$\elr$}\rput(.5,.5){$\eur$}
\rput(-.5,-.5){$\ev$}\rput(-.5,.5){$\eul$}
\rput(0,-1.5){at $z= w$}\rput(.2,-1.88){$\Gamma_+=\Gamma_-$}
}
\end{pspicture}\]
yield for $A'=A_z$ the exact sequences (\ref{seq2:hatbarz}), (\ref{seq2:hatV}) and (\ref{seq2:hatbarw}).
\end{proof}                          

An analysis for $\nabla_{A_z}^*\nabla_{A_z}$ brings a very similar
wall crossing formula
\[\ind(\nabla_{A_z}^*\nabla_{A_z},\eul)-\ind(\nabla_{A_z}^*\nabla_{A_z},\elr)
=\begin{cases} 2\dim \Vlam_0,&\text{ for }\Gamma_+\neq \Gamma_-;\\
   4\dim \Vlam_0,&\text{ for }\Gamma_+= \Gamma_-.\end{cases}\]
However, since $\nabla_{A_z}^*\nabla_{A_z}$ is
self-adjoint, $\ind(\nabla_{A_z}^*\nabla_{A_z},\eul)=
-\ind(\nabla_{A_z}^*\nabla_{A_z},\elr)$, whence
\[\rk\HHH=
\begin{cases} \dim \Vlam_0,&\text{ for }\Gamma_+\neq \Gamma_-;\\
   2\dim \Vlam_0,&\text{ for }\Gamma_+= \Gamma_-.\end{cases}\]
Using Equation (\ref{rkV}), we can even say
\[\rk\HHH=
\begin{cases}
  2,&\text{ for }\Gamma_+\neq \Gamma_-\text{ and }2w\not\in\Lambda^*;\\
  4,&\text{ for }\Gamma_+\neq \Gamma_-\text{ and }2w\in\Lambda^*,
     \text{ or }\Gamma_+=\Gamma_-\text{ and }2w\not\in\Lambda^*;\\
  8,&\text{ for }\Gamma_+= \Gamma_-\text{ and }2w\in\Lambda^*.
\end{cases}\]

Similarly, we have for the Laplacian the following isomorphisms:
\begin{align}
0 \longrightarrow  \HHH_z \longrightarrow \Vlam_\lambda\oplus
  \Vlam_{-\lambda}\longrightarrow 0,
   \quad\quad&\text{ for }z\neq w\text{ and when }\Gamma_+\neq \Gamma_-,
   \label{seq:HVz}\\
0\longrightarrow \HHH_w\longrightarrow  \Vlam_0\longrightarrow 0,\quad\quad&
\text{ when }\Gamma_+\neq \Gamma_-,\label{seq:HV0}\\
0 \longrightarrow  \HHH_z \longrightarrow \bigl(\Vlam_\lambda\oplus
  \Vlam_{-\lambda}\bigr)^{\oplus2}\longrightarrow 0,\quad\quad
   & \text{ for }z\neq w\text{ and when }\Gamma_+= \Gamma_-,
   \label{seq2:HVz}\\
0\longrightarrow \HHH_w \longrightarrow  \Vlam_0\oplus \Vlam_0\longrightarrow 0,
\quad\quad&\text{ when }\Gamma_+= \Gamma_-.\label{seq2:HV0}
\end{align}

Bringing all of those sequences together allows us to conclude the following.

\begin{theorem}On $\Bw$, we have     
\[\bE=\vE\oplus\DD\HHH.\]
\end{theorem}

\begin{proof}Denote $\Wlam_\lambda$ the space  
\[\Wlam_\lambda:=\begin{cases} \Vlam_\lambda\oplus \Vlam_{-\lambda},&\text{ if }\Gamma_+=\Gamma_-;\\ \Vlam_\lambda,&\text{ if }\Gamma_+\neq\Gamma_-.\end{cases}\]
Let $p\colon \Wlam_\lambda\oplus \Wlam_{-\lambda}\to \Wlam_\lambda$ denote the map
$p(a,b)=2\lambda a$.

For $\lambda\neq 0$, we use the Snake Lemma\index{snake lemma} on the diagram
\[\begin{CD}&&0@>>> \HHH @>>> \Wlam_\lambda\oplus  \Wlam_{-\lambda} @>>> &0\\
&&@VVV @V{\DD}VV @VV{p}V\\
0@>>> \Ehat @>>>  \bE @>>> \Wlam_\lambda
@>>> 0
\end{CD}\]
coming from Sequences (\ref{seq:hatbarz}), (\ref{seq2:hatbarz}), (\ref{seq:HVz}), and (\ref{seq2:HVz}), to produce an exact sequence
\begin{alignat}6
&\ker(0)&\longrightarrow
&\ker(\DD)&\longrightarrow
&\ker(p)&\longrightarrow
&\coker(0)&\longrightarrow
&\coker(\DD)&\longrightarrow
&\coker(p)\nonumber\\
&\phantom{ker} 0&\longrightarrow
&\phantom{ke} \bK_z&\longrightarrow
&\phantom{k} \Wlam_{-\lambda}&\longrightarrow
& \phantom{cok}\Ehat_z&\longrightarrow
&\coker(\DD)&\longrightarrow
&\phantom{co} 0\label{seq:fromsnake}
\end{alignat}
Note that the map $\Ehat\to\coker(\DD)$ being surjective forces $\bE$
to be spanned by $\Ehat$ and $\DD\HHH$.

Sequences (\ref{seq:hatV}) and (\ref{seq2:hatV}) imply
\[\dim \Ehat_z = \dim\vE_z +\dim \Wlam_\lambda -\dim\bK_z\]
while Sequences (\ref{seq:hatbarz}) and (\ref{seq2:hatbarz}) imply
\[\dim \bE_z=\dim \Ehat_z+\dim \Wlam_\lambda.\]
Thus
\[\dim \bE_z=\dim \vE_z+2\dim \Wlam_\lambda-\dim \bK_z=\dim \vE_z+\dim \DD\HHH.\]
Since Lemma \ref{lemma:kill} guarantees that $\scp{\DD\HHH,\vE}=\{0\}$, we have
$\Ehat\cap\DD\HHH$ perpendicular to $\vE$ for the $L^2$ inner product.  Hence
$\DD\HHH\cap\vE=\{0\}$, and $\bE_z=\vE_z\oplus\DD\HHH$.

It remains to prove the theorem for $z=w$.  We already know $\vE_w=\Ehat_w$
and $\DD\HHH_w\subset\bE_w$.  We also know from Sequences (\ref{seq:hatbarw})
and (\ref{seq2:hatbarw}) that
\begin{align*}\dim\bE_w&=\dim\Ehat_w+\dim \Wlam_0-\dim \bK_w\\
&=\dim\vE_w +\dim\DD\HHH_w.\end{align*}
We therefore only have to prove that the intersection
$\Ehat_w\cap\DD_{A_w}\HHH_w$ is $\{0\}$ to complete the proof.

The asymptotic behavior of $\phi\in\HHH_w$ is
\[\phi=\begin{cases} t\phi_0^++\phi_1^++o(1), &\text{ as }t\to\infty;\\
                     t\phi_0^-+\phi_1^-+o(1), &\text{ as }t\to-\infty;
  \end{cases}\]
for some $\phi_0^\pm,\phi_1^\pm\in \Vlam_0$.  
If $\Gamma_+\neq\Gamma_-$,
we must have $\phi_0^-=\phi_1^-=0$, as $w$ is associated to $\Gamma_+$.

The asymptotic behavior of $\DD_{A_w}\phi$ is
\[\DD_{A_w}\phi=\begin{cases} \phi_0^++o(1), &\text{ as }t\to\infty;\\
                     \phi_0^-+o(1), &\text{ as }t\to-\infty.
  \end{cases}\]
Suppose $\DD_{A_w}\phi\in L^2$.  Then
\begin{align*} \n{\DD_{A_w}\phi}{L^2}^2 &=\scp{\DD_{A_w}^*\DD_{A_w}\phi,\phi}
  +\lim_{t\to\infty}\scp{\DD_{A_w}\phi,\cl(\partial_t)\phi}
  +\lim_{t\to-\infty}\scp{\DD_{A_w}\phi,\cl(\partial_t)\phi}\\
&=\scp{\phi^+_0,\phi^+_1}+\lim_{t\to\infty}t|\phi^+_0|^2
  -\scp{\phi^-_0,\phi^-_1}-\lim_{t\to-\infty}t|\phi^-_0|^2.
\end{align*}
For $\n{\DD_{A_w}\phi}{L^2}$ to be finite, we must get rid of the limits, thus
forcing $\phi_0^\pm=0$ and consequently we have $\DD_{A_w}\phi=0$. \cqfd
\end{proof}

For a continuous family of Fredholm operators, like $\DD_{A_z}$ on $L^2_\eul$
parameterized on $\Bw$, the dimension of the kernel can only drop in a small
neighborhood of a given point, it cannot increase.  However, not any random
behavior is acceptable.

\begin{lemma}[see \citep{p. 241}{kato}] \label{lemma:kato} Let
$T\colon X\to Y$ be Fredholm
and $S\colon X\to Y$ a bounded operator.  Then
the operator $T+tS$ is Fredholm and $\dim\ker(T+tS)$ is constant
for small $|t|>0$.\end{lemma}

We obviously use this lemma
with $T=\DD_{A_w}, X=\Wotd{\eul},Y=L^2_{\eul}$, and $S=\cl(e)$ for some
direction $e\in\R^3$. Let's note that three scenarios are possible.
\begin{enumerate}\label{scenarios}
\item $\dim\bK_z$ is constant on a neighborhood around $w$, say $\Bw$;
   \label{case:constant}
\item $\dim\bK_z$ is constant for $z\in\Bw\setminus\{w\}$, but is smaller
than $\dim\bK_w$;
   \label{case:alldrop}
\item $\dim\bK_{w+\lambda e}\neq\dim\bK_{w+\lambda'e'}$ for
small $\lambda,\lambda' >0$ and some $e\neq e'$.
   \label{case:droponeray}
\end{enumerate}

\section{Asymptotic of the Higgs field} 
\label{sec:Higgsfield}                  
We now study the behavior of the Higgs field $\Higgs$ as
$z$ approaches of a given element $w$  of $W$.
We know $w$ is associated to the limit $\Gamma$ of $A$ at $\infty$ or
$-\infty$, in the sense that $\Gamma$ splits $E$ as $L_w\oplus L_{-w}$.
Without loss of generality, we suppose
\[\Gamma_+=\Gamma.\]

When $\Gamma_+\neq \Gamma_-$, and for $2\pi |z-w|<\epsilon$, notice that
\begin{align*}
\bE_z &= L^2_\eb\cap\ker(\DD_{A_z}^*) = L^2_\ep\cap\ker(\DD_{A_z}^*)
      = L^2_\eur\cap\ker(\DD_{A_z}^*), \text{ and}\\
\vE_z &= L^2_\ev\cap\ker(\DD_{A_z}^*) = L^2_\emo\cap\ker(\DD_{A_z^*})
       = L^2_\elr\cap\ker(\DD_{A_z}^*).\end{align*}
When $\Gamma_+= \Gamma_-$, those spaces are a priori all different.

\begin{theorem}\label{thm:constantrank} 
On a closed ball
$\Bw$ around $w$,
there exists families of operators
$\Pb$ and $\Pv$, bounded independently of
$z$ , such that
\begin{equation}\label{HiggsSplitting}
      \Higgs=\frac {-i}{2|z-w|}\Pb+\Pv.\end{equation}
Furthermore, $\Pb$ is the $L^2$-orthogonal projection on
$\DD_{A_z}\HHH_z\cap \Ehat_z$.
\end{theorem}

\begin{proof}    
Obviously, $\vE$ supports many different
norms, and amongst those are the $L^2$ and $L^2_\elr$ norms.
For $\phi\in\vE_z$,  observe that
\[\n{t\phi}{L^2}\leq C_\epsilon\n\phi{L^2_\elr}.\]
We would really like to bound this last quantity by a multiple of
$\n\phi{L^2}$.

Let $Q$ denote the projection $L^2\to \Eh_w$.  Of course,
since $L^2_\elr\subset L^2$, the projection is also defined on $L^2_\elr$.
Let $\Eh_w^\perp$ be the $L^2$-orthogonal complement, and
$\Eh_w^0=\Eh_w^\perp\cap L^2_\elr$.  In fact, we have
\[ L^2_\elr = \Eh_w\oplus \Eh_w^0\]
since at $w$, we have $\Eh_w=\vE_w$.

Since $\DD_{A_w}^*$ is injective on $\Eh_w^0$, there is a constant such
that
\[\n{u}{L^2_\elr}\leq C\n{\DD_{A_w}^*u}{L^2_\elr}\text{ for }u\in \Eh_w^0.\]
But then for $u\in \vE_{w+\lambda e}$, we have
\begin{align*}
 \n{u}{L^2_\elr} &\leq \n{Qu}{L^2_\elr}+\n{(1-Q)u}{L^2_\elr}\\
   &\leq \n{Qu}{L^2_\elr}+C\n{\DD_{A_w}^*(1-Q)u}{L^2_\elr}\\
   &= \n{Qu}{L^2_\elr}+C\n{\DD_{A_w}^*u}{L^2_\elr}\\
   &= \n{Qu}{L^2_\elr}+C\lambda\n{u}{L^2_\elr}.
\end{align*}
Hence for $\lambda$ small enough,
\[\n{u}{L^2_\elr}\leq 2\n{Qu}{L^2_\elr}.\]

Of course, since $\Eh_w$ is finite dimensional, there exists a constant $C$ for which
$\n{Qu}{L^2_\elr}
\leq C\n{Qu}{L^2}$ and thus for $u\in\vE_z$ with $z$ close to $w$,
\[\n{tu}{L^2}\leq C_\epsilon\n{u}{L^2_\elr}\leq 2C_\epsilon \n{Qu}{L^2_\elr}
\leq C\n{Qu}{L^2}\leq  C\n{u}{L^2}.\]

Denote $\Prv$ the $L^2$-orthogonal projection of $\Ehat$
on $\vE$. We just proved that
\[\Higgs\circ\Prv\text{ is bounded independently of }z\in\Bw.\]
It is part of the map $\Pv$ announced in the statement of the theorem.

One of the crucial feature of this proof is our ability to find a uniform
bound for $m_t$ on $\vE$.

As suggested above, let $\Pp$ denote the $L^2$-orthogonal projection on
$\DD_{A_z}\HHH_z\cap \Ehat_z$.  Then
\begin{align*}
\Higgs=-2\pi i  P m_t & = \Higgs\Prv-2\pi i\bigl(\Prv+\Pp\bigr) m_t \Pp\\
   &=\Higgs\Prv+2\pi i\Prv m_t \Pp-2\pi i\Pp m_t \Pp.\end{align*}

For $\phi_1\in \vE$, and $\phi_2\in V$, we have
$\scp{\phi_1,t\Pp\phi_2}=\scp{t\phi_1,\Pp\phi_2}$.
Thus $\Prv m_t \Pp$ is also bounded independently of $z\in\Bw$.


It remains only to analyze $\Pp m_t \Pp$.
Pick a vector $e\in\R^3$ of length $1$.   Let
\[{\cal R}=\{w+\frac\lambda{2\pi} e\}\subset \Bw\]
be a ray inside $\Bw$ emerging from $w$.  As the notation suggests, we
parameterize this ray by $\lambda=2\pi |z-w|$.  Pick a
family $\phi_z\in\DD_{A_z}\HHH_z$ for $z\in {\cal R}$, with
\begin{gather} \phi_z\in \Ehat_z \text{ for }\lambda>0,\notag\\
  \n{\phi_z}{L^2_\eul}=1.\label{normalizationphi}\end{gather}

But then,
\[\n{\phi_z}{L^2}\to\infty\text{ as }\lambda\to 0.\]

To prove this claim, suppose it is not true.  Then there is a subsequence
   $\phi_{z_j}\rightharpoonup\tilde\phi_w$ weakly in $L^2$.  Hence
$\scp{\phi_{z_j},f}\to\scp{\tilde\phi_w,f}$ for all $f\in L^2$,
in particular for all $f\in L_{\elr}^2={(L_\eul^2)}^*$,
whence $\phi_{z_j}\rightharpoonup\tilde\phi_w$ weakly in
$L_\eul^2$.  Since $\phi_z\to \phi_w$ in $L^2_\eul$,
we have $\tilde\phi_w=\phi_w$, which is impossible as
$\tilde\phi_w$ is in $L^2$ while $\phi_w$ is not.

Because $\Gamma_w$ is independent of $t$, and because $-\epsilon$ is not
an eigenvalue of $D_{\Gamma_w}$,
the operator $\DD^*_{\Gamma_w}$ is an isomorphism $\Wotd\ev\to L^2_\ev$, and
$\Wotd\eur\to L^2_\eur$, hence
there exist a constant $C$ such that
\begin{align}
    \n{u}{\Wotd\ev}\leq C\n{\DD_{\Gamma_w}^*u}{L^2_\ev},&\quad\text{ for }u
            \in \Wotd\ev,\label{Isomcheck}\\
    \n{u}{\Wotd\eur}\leq C\n{\DD_{\Gamma_w}^*u}{L^2_\eur},&\quad\text{ for }u
             \in \Wotd\eur.\label{Isomur}
\end{align}

Because $\phi_z\in\Ehat_z$ for $\lambda>0$,
for $t>0$, we can write
$\phi_z=e^{-\lambda t}\psi_{-\lambda} + g_z$ for some
eigenvector $\psi_{-\lambda}$ of eigenvalue $-\lambda$ of
$D_{\Gamma_z}$
and some $g_z\in \Wotd{-\epsilon}(\Tp)$.
When $\Gamma_-=\Gamma_+$, and for $t<0$,
we can write
$\phi_z=e^{\lambda t}\psi_{\lambda} + j_z$ for some
eigenvector $\psi_{\lambda}$ of eigenvalue $\lambda$ of
$D_{\Gamma_z}$
and some $j_z\in \Wotd{\epsilon}(\Tm)$.

While $g_z$ and $j_z$ appear to be defined only for $t>0$ and $t<0$
respectively,  let's define them globally on $\RT$ by
$g_z=\phi_z-e^{-\lambda t}\psi_{-\lambda}$ and $j_z=\phi_z-e^{\lambda t}
\psi_{\lambda}$.

Notice that
\begin{equation}\label{eqn:Dgg}
\DD_{\Gamma_z}^*g_z = \DD_{\Gamma_z}^*\phi_z
  =(\DD_{\Gamma_z}^*-\DD_{A_z}^*)\phi_z=\cl(\Gamma-A)\phi_z,
\end{equation}
and similarly
\begin{equation}\label{eqn:Dgj}
\DD_{\Gamma_z}^*j_z = \DD_{\Gamma_z}^*\phi_z
  =(\DD_{\Gamma_z}^*-\DD_{A_z}^*)\phi_z=\cl(\Gamma-A)\phi_z,
\end{equation}

Overall, there is a constant such that
$|\cl(A-\Gamma)|\leq C\sigma_{(0,\beta)}$,
and this estimate can be improved to
$|\cl(A-\Gamma)|\leq C\sigma_{(-\beta,\beta)}$ when $\Gamma_-=\Gamma_+$.
Hence $\cl(A-\Gamma)$ gives a bounded map $L^2_\eb\to L^2_\ev$ in all cases
and $L^2_\eb\to L^2_\eur$ when $\Gamma_-=\Gamma_+$.
Thus Equation (\ref{eqn:Dgg})  yields
\begin{equation}\label{boundclleft}
     \n{\DD_{\Gamma_z}^*g_z}{L^2_\ev}\leq C\n{\phi_z}{L^2_\eb},
\end{equation}
and for the special case $\Gamma_-=\Gamma_+$, Equation (\ref{eqn:Dgj}) yields
\begin{equation}\label{boundcltop}
     \n{\DD_{\Gamma_z}^*j_z}{L^2_\eur}\leq C\n{\phi_z}{L^2_\eb}.
\end{equation}

{}From Equations (\ref{Isomcheck}), and (\ref{boundclleft}),
we derive
\begin{align*}  \n{g_z}{\Wotd\ev} &\leq  C\n{\DD_{\Gamma_w}^*g_z}{L^2_\ev}\\
&=C\n{\DD_{\Gamma_z}^*g_z+\lambda \cl(e) g_z}{L^2_\ev}\\
&\leq C\n{\phi_z}{L^2_\eb}+C\lambda\n{g_z}{L^2_\ev},
\end{align*}
After rearranging, we notice that
  $\n{g_z}{\Wotd\ev}$ is bounded independently of small $z$,
and
similarly $\n{j_z}{\Wotd\eur}$ is bounded independently of small z.
This last fact is also true for $\Gamma_-\neq\Gamma_+$, for in that
case $j_z=\phi_z$ and its $L^2_\eb$-norm is equivalent to the $L^2_\eur$-norm,
as both as defined on $\bE$ over $\Bw$.

While it is agreeable to work with a smooth splitting, nothing prevents us
{}from considering the  functions
\[ h_\lambda=\begin{cases}
     e^{\lambda t}\psi_\lambda, &\text{ for }t<0,\\
     e^{-\lambda t}\psi_{-\lambda},  &\text{ for }t>0,\end{cases}
\quad\quad\text{ and }\quad\quad r_z = \begin{cases} j_z, &\text{ for }t<0,\\
                           g_z,  &\text{ for }t>0,\end{cases}
\]
and the associate splitting
\[\phi_z = h_\lambda+ r_z.\]

That
$\n{r_z}{L^2_\elr}$ is bounded independently of small z.
follows from the similar fact concerning $g_z$ and $j_z$.

Consider the families
\begin{align*} \bar\phi_z&:={\phi_z}/{\n{\phi_z}{L^2}},\\
\bar h_\lambda&:=h_\lambda/  \n{\phi_z}{L^2},\\
\bar r_z&:=r_z/  \n{\phi_z}{L^2}.\end{align*}

Since $\n{\phi_z}{L^2}\to\infty$ and $\n{r_z}{L^2_\elr}$ is bounded, we have
$\n{\bar r_z}{L^2_\elr}\to 0$  as $\lambda\to 0$, and
a fortiori, $\n{\bar r_z}{L^2}\to 0$.
The triangle inequality then guarantees
\[\bigl|\n{\bar h_\lambda}{L^2}-\n{\bar r_z}{L^2}\bigr|
   \leq \n{\bar \phi_z}{L^2}\leq
  \n{\bar h_\lambda}{L^2}+\n{\bar r_z}{L^2}.\]
Since $\n{\bar \phi_z}{L^2}=1$, and $\n{\bar r_z}{L^2}\to 0$,
we must have
\[\n{\bar h_\lambda}{L^2}\to 1 \text{ as }\lambda\to 0.\]

Let's now come back to our main worry.  We study
\[\scp{t\bar\phi_z,\bar\phi_z}
=\scp{t\bar h_\lambda,\bar h_\lambda}+2\scp{\bar h_\lambda,t\bar r_z}+
  \scp{t\bar r_z,\bar r_z}.\]
The last two terms are bounded by a multiple of $\n{t\bar r_z}{L^2}$.
But
\[\n{t\bar r_z}{L^2}\leq C \n{\bar r_z}{L^2_\elr},\]
hence it is going to $0$.

As for the first term, we have
\begin{align*}
\scp{t\bar h_\lambda,\bar h_\lambda}
   &=\frac1{\n{\phi_\lambda}{L^2}^2}
     \Bigl(\int_0^\infty te^{-2\lambda t}|\psi_{-\lambda}|^2
   +\int_{-\infty}^0 te^{2\lambda t}|\psi_{\lambda}|^2\Bigr)\\
   &=\frac1{2\lambda}\frac1{\n{\phi_\lambda}{L^2}^2}
     \Bigl(\int_0^\infty e^{-2\lambda t}|\psi_{-\lambda}|^2+
       \int_{-\infty}^0 e^{2\lambda t}|\psi_{\lambda}|^2\Bigr)\\
   &=\frac1{2\lambda}\n{\bar h_\lambda}{L^2}^2,\end{align*}
hence
\[\scp{t\bar\phi_\lambda,\bar\phi_\lambda}=\frac1{2\lambda}+o(1)
                                             \text{ as }\lambda\to 0.\]

Suppose now $\bar\phi^1_z$ and $\bar\phi^2_z$ are two such families, but so
that \[\scp{\bar\phi^1_z,\bar\phi^2_z}_{L^2}=0.\]
Then
\begin{align*}\scp{t\bar\phi^1_z,\bar\phi^2_z}
   &=\scp{t\bar h^1_\lambda,\bar h^2_\lambda}
     +\scp{\bar h^1_\lambda,t\bar r^2_z}+\scp{t\bar r^1_z,\bar h^2_\lambda}
+\scp{t\bar r^1_z,\bar r^2_z}\\
&=\frac1{2\lambda}\scp{\bar h^1_\lambda,\bar h^2_\lambda} +o(1),
\end{align*}
and of course $\scp{\bar h^1_\lambda,\bar h^2_\lambda}\to 0$, hence
the result.
\end{proof}

Finally, let's note that in fact, Scenario
\ref{case:droponeray} of page~\pageref{scenarios} cannot happen.
We can take the trace of $(\Ah,\Higgs)$ to obtain an abelian
monopole $(b,\varphi)$ on $\Bw\setminus\{w\}$.    The Bogomolny equation
reduces to
\[d\, \varphi=*d\, b,\]
and thus $\Delta\varphi=0$.  Since $\varphi$ is harmonic, not every possible
behavior  is acceptable as $z\to w$.  For one thing, there is a unique
set of homogeneous harmonic polynomials $p_m$ and $q_m$ of degree $m$
which give a decomposition of $\varphi$ on $\Bw\setminus\{w\}$ as a Laurent
series
\[\varphi=\sum_{m=0}^\infty p_m(z-w)+\sum_{m=0}^\infty
     \frac{q_m(z-w)}{|z-w|^{2m+1}};\]
see for example \cite[Thm 10.1, p. 209]{axler}.

Whether or not the rank is  constant, we can find for any sequence of points
approaching $w$ a subsequence
of points $z_j\to w$ for which the decomposition of Equation
\ref{HiggsSplitting} is valid.  We then have
\[\lim_{j\to\infty}2|z_j-w|\varphi_{z_j}=i\dim\DD_{A_{z_j}}\HHH_{z_j}
=i(\rk \HHH - \dim \bK_{z_j}).\]
By the Laurent series decomposition given above, this number
must be the same in any way we approach $w$, hence $\dim\bK_z$ must
be constant on $\Bw\setminus\{w\}$.


\newcommand{\noopsort}[1]{} \def\cprime{$'$}
\providecommand{\bysame}{\leavevmode\hbox to3em{\hrulefill}\thinspace}
\providecommand{\href}[2]{#2}

{\sc Department of Mathematics\\
McGill University\\
805 Sherbrooke West\\
Montr\'eal, QC H3A 2K6\\
Canada}\\
\emph{E-mail address:} \texttt{benoit@alum.mit.edu}\\

{\sc Received by Communications in Analysis and Geometry\\ on October 20, 2004.}

\end{document}